\documentclass[10 pt]{amsart}

\usepackage{latexsym, amsmath, amssymb, longtable, booktabs,amscd,microtype,booktabs,cases}
\usepackage[square]{natbib}
\usepackage{hyperref}
\bibpunct{[}{]}{,}{n}{}{;} 

\numberwithin{equation}{section}
\begin{document}

\newcommand\A{\mathbb{A}}
\newcommand\C{\mathbb{C}}
\newcommand\G{\mathbb{G}}
\newcommand\N{\mathbb{N}}
\newcommand\T{\mathbb{T}}
\newcommand\sO{\mathcal{O}}
\newcommand\sE{{\mathcal{E}}}
\newcommand\tE{{\mathbb{E}}}
\newcommand\sF{{\mathcal{F}}}
\newcommand\sG{{\mathcal{G}}}
\newcommand\GL{{\mathrm{GL}}}
\newcommand\HH{{\mathrm H}}
\newcommand\mM{{\mathrm M}}
\newcommand\fS{\mathfrak{S}}
\newcommand\fP{\mathfrak{P}}
\newcommand\fQ{\mathfrak{Q}}
\newcommand\Qbar{{\bar{\Q}}}
\newcommand\sQ{{\mathcal{Q}}}
\newcommand\sP{{\mathbb{P}}}
\newcommand{\Q}{\mathbb{Q}}
\newcommand{\tH}{\mathbb{H}}
\newcommand{\Z}{\mathbb{Z}}
\newcommand{\R}{\mathbb{R}}
\newcommand{\F}{\mathbb{F}}
\newcommand\gP{\mathfrak{P}}
\newcommand\Gal{{\mathrm {Gal}}}
\newcommand\SL{{\mathrm {SL}}}
\newcommand\Hom{{\mathrm {Hom}}}
\newcommand{\legendre}[2] {\left(\frac{#1}{#2}\right)}
\newcommand\iso{{\> \simeq \>}}

\newtheorem{thm}{Theorem}
\newtheorem{theorem}[thm]{Theorem}
\newtheorem{cor}[thm]{Corollary}
\newtheorem{conj}[thm]{Conjecture}
\newtheorem{prop}[thm]{Proposition}
\newtheorem{lemma}[thm]{Lemma}
\theoremstyle{definition}
\newtheorem{definition}[thm]{Definition}
\theoremstyle{remark}
\newtheorem{remark}[thm]{Remark}
\newtheorem{example}[thm]{Example}
\newtheorem{claim}[thm]{Claim}

\newtheorem{lem}[thm]{Lemma}

\theoremstyle{definition}
\newtheorem{dfn}{Definition}

\theoremstyle{remark}

\theoremstyle{remark}
\newtheorem*{fact}{Fact}
% type user-defined commands here
\makeatletter
\def\imod#1{\allowbreak\mkern10mu({\operator@font mod}\,\,#1)}
\makeatother

\title{The Eisenstein elements of modular symbols for level product of two distinct odd primes}
\author{Debargha Banerjee}
\address{INDIAN INSTITUTE OF SCIENCE EDUCATION AND RESEARCH, PUNE, INDIA}
\author{SriLakshmi Krishnamoorty}
\address{INDIAN INSTITUTE OF TECHNOLOGY, CHENNAI, INDIA}
\begin{abstract}
We explicitly write down the {\it Eisenstein elements} inside the space of modular symbols for  
Eisenstein series with integer coefficients for the congruence subgroups $\Gamma_0(pq)$ with $p$ and $q$ distinct odd primes, giving an answer 
to a question of Merel in these cases. 
We also compute the winding elements explicitly 
for these congruence subgroups. Our results are explicit versions of the Manin-Drinfeld Theorem
[Thm.\,\ref{ManinDrinfeld}].
\end{abstract}

\maketitle
 
\section{Introduction}
\label{Introduction}
In his landmark paper on Eisenstein ideals, Mazur  studied torsion points of elliptic curves over $\Q$ and gave a list 
of possible torsion subgroups of elliptic curves [cf. Thm. 8, \cite{MR488287}]. 
In \cite{MR1405312}, Merel 
wrote down modular symbols for the congruence subgroups $\Gamma_0(p)$ for any odd prime $p$ 
that correspond to differential forms of third kind on the modular curves.
He then use these modular symbols to give an uniform upper bound of the torsion 
points of elliptic curves over any number fields in terms of extension degrees of these number fields
 \cite{MR1369424}. 
The explicit expressions of winding elements for prime level of \cite{MR1405312} are being used by Calegari and Emerton to 
study the ramifications of Hecke algebras at the Eisenstein primes \cite{MR2129709}. 
 Several authors afterwards
studied the torsion points of elliptic curves over number fields using modular symbols.

In the present paper, we  study elements of relative homology groups of the 
modular curve $X_0(pq)$  that correspond to differential forms of 
third kind with $p$ and $q$ distinct odd primes.  
As a consequence, we give an ``effective" proof 
of the Manin-Drinfeld theorem [Thm.\,\ref{ManinDrinfeld}] for the special case of the image in $\HH_1(X_0(pq), \R)$ of the path in $\HH_1(X_0(pq), \partial(X_0(pq)), \Z)$ joining $0$  and $i\infty$. Since the algebraic part of the special values 
of $L$-function are obtained by integrating differential forms on these modular symbols, our explicit expression of the winding elements  should be 
useful to understand the algebraic parts of the special values at $1$ of the $L$-functions of the quotient Jacobian of modular curves for the congruence 
subgroup $\Gamma_0(pq)$ \cite{amod}.  

For $N \in \{p,q,pq\}$,  consider the basis $E_N$ of $E_2(\Gamma_0(pq))$ [\S\,\ref{Eisensteinseries}] for which all the Fourier coefficients at $i\infty$ belong to $\Z$.
 The meromorphic differential forms $E_N(z)dz$ are of third kind on the Riemann surface
$X_0(pq)$ but of first kind on the non-compact Riemann surface $Y_0(pq)$.

Let $\xi : \SL_2(\Z) \rightarrow \HH_1(X_0(pq),\text{cusps},\Z)$ 
be the Manin map [\S~\ref{Mod}]. 
 For any two coprime integers $u$ and $v$ with $v \geq 1$, let $S(u,v) \in \Z$ be the Dedekind sum
 [cf \ref{Peiods}].  If 
$g \in \sP^1(\Z \slash pq \Z)$ is not of the form 
 $(\pm 1 ,1),  ( \pm1  \pm kx, 1)$ or  $(1, \pm 1 \pm kx) $ 
with $x$ one of the prime $p$ or $q$, then we can write it as $(r-1,r+1)$.

Let $\delta_r$ be $1$ or $0$ depending on $r$ is odd or even. 
 For any integer $k$, let $s_k=(k+(\delta_k-1)pq)$ be an odd integer.
 Choose integers 
$s, s'$ and $l, l'$ such that $l(s_kx+2)-2spq=1$ and $l's_kx-2s'\frac{pq}{x}=1$.  Let $\gamma_1^{x,k}=\left(\begin{smallmatrix}
1+4spq &  -2l \\
- 4s(s_kx+2)pq& 1+4spq\\
\end{smallmatrix}\right)$ 
and $\gamma^{x,k}_2=\left(\begin{smallmatrix}
1+4s'\frac{pq}{x} &  -2l' \\
- 4s'(s_k)pq& 1+4s'\frac{pq}{x}\\
\end{smallmatrix}\right)$ be two matrices [cf.~\ Lemma \ref{divisibleinteger}]. For $l=1,2$, consider the integers 
\[
P_N(\gamma^{x,k}_{l})=sgn(t(\gamma^{x,k}_{l}))[2 (S(s(\gamma^{x,k}_l),|t(\gamma^{x,k}_l)|N)-S(s(\gamma^{x,k}_l), |t(\gamma^{x,k}_l)|))
\]
\[
-S(s(\gamma^{x,k}_l),\frac{|t(\gamma^{x,k}_l|)}{2}N)+S(s(\gamma^{x,k}_l), \frac{|t(\gamma^{x,k}_l)|}{2})]
\]
with 
\[
 s(\gamma_1^{x,k})=1-4spq(1+s_kx), t(\gamma_1^{x,k})=-2(l-2s(s_kx+2)pq)
 \]
  and 
 \[
 s(\gamma_2^{x,k})=1-4s'pq(s_k-\frac{1}{x}), t(\gamma_2^{x,k})=-2(l'-2s's_kpq).
 \]
Define the function $F_N: \sP^1(\Z \slash pq \Z) \rightarrow \Z$ by  
\begin{equation*}
F_{N}(g) =
\begin{cases}
2(S(r,N)-2S(r,2N)) & \text{if $g=(r-1,r+1)$,} \\
P_N(\gamma_1^{x,k})-P_N(\gamma^{x,k}_2) &  \text{if $g=(1+kx,1)$ or $g=(-1-kx,1)$,} \\
-P_N(\gamma_1^{x,k})+P_N(\gamma^{x,k}_2) &  \text{if $g=(1,1+kx)$ or $g=(1,-1-kx)$,} \\
0 &  \text{if $g=(\pm 1,1)$.} \\
\end{cases}
\end{equation*}
\begin{thm}
\label{Main-thm}
The modular symbol
$$\sE_{E_N}= \sum_{g \in \sP^1(\Z \slash pq \Z)} F_{E_N}(g) \xi(g)$$
in $\HH_1(X_0(pq), \partial(X_0(pq), \Z)$ is the Eisenstein element
 [\S ~\ref{Eisensteinelements}] corresponding to the Eisenstein series $E_N \in E_2(\Gamma_0(pq))$.
\end{thm}
In \cite{debargha}, a description is given of Eisenstein elements in terms of certain integrals 
for $M=p^2$. In this article, we give an explicit description in terms of two matrices $\gamma_1^{x,k}$ and 
$\gamma_2^{x,k}$.
Let $\overline{B_1} :\R \rightarrow \R$ be the periodic first Bernoulli polynomial. For the Eisenstein series 
$E_{pq}$ [\S\,\,\ref{Eisensteinseries}], we write down the Eisenstein elements more explicitly if 
$g=(r-1,r+1)$. Replacing $p$ with $pq$ [Lemma 4,~\cite{MR1405312}], we write 
\[
F_{pq}((r-1,r+1))=\sum_{h=0}^{pq-1} \overline{B_1} (\frac{hr}{2pq}).
\]
Recall the concept of the {\it winding elements} [\S ~\ref{Winding}]. We write down the 
explicit expression of the winding elements for the congruence subgroup $\Gamma_0(pq)$.  
\begin{cor}
\label{cor}
 $$(1-pq) e_{pq}=\sum_{x \in(\Z/pq\Z)^*} F_{pq}((1,x)) \{0,\frac{1}{x}\}.$$
\end{cor}
Note that if $\nu=gcd(pq-1,12)$ and $n=\frac{pq-1}{\nu}$, then a multiple of winding element $ne_{pq}$ belongs to 
$\HH_1(X_0(pq), \Z)$.
Manin-Drinfeld proved that the modular symbol 
$\{0, \infty \} \in \HH_1(X_0(N), \Q)$ using the theory of suitable Hecke operators 
acting on modular curve $X_0(N)/ \Q$.
In this paper, we follow the approach of Merel [cf.~\cite{MR1405312}, Prop. 11].
 Our explicit expression of winding element 
should be useful to understand the algebraic part of the {\it special values
of $L$-functions} [cf.~\cite{amod}, p. 26].

Since Hecke operators are defined over $\Q$, there is a possibility that we can find the Eisenstein elements for the congruence 
subgroups of odd level in a completely different method without using boundary computations. 
It is tempting to remark that our method should generalize to the congruence subgroup $\Gamma_0(N)$ 
atleast if $N$ is squarefree and odd. Unfortunately,  generalizing our method is equivalent 
to explicit understanding  of boundary homologies  of modular curves defined over rationals.  For instance, if $N=pqr$ with $p, q, r$ three distinct primes then there are $8$ cusps. 
Since in these cases, there are more cusps the computations of boundaries become much more tedious. 
One of the author wish to tackle the difficulty using the ``level" of the cusps in a future article. 
\section{Acknowledgements}
This paper owes it's existence to several e-mail communications, encouragements and several valuable suggestions of 
Professor Lo\"ic Merel. We would like to thank IMSc, Chennai for providing excellent working conditions. The second author would like to thank MPIM, Germany
for the great hospitality during her visit. She was supported by a DST-INSPIRE grant. We wish to thank the anonymous  referee for 
several useful comments that helped us to improve the mathematical content of this article. 
\section{Modular Symbols}
\label{Mod}
Let $\tH \cup \sP^1(\Q)=\overline{\tH}$ and $\Gamma \subset \SL_2(\Z)$ be a congruence subgroup. The topological space 
$X_{\Gamma}(\C) =\Gamma \backslash \overline{\tH}$ has a natural structure of a smooth compact  Riemann surface and 
consider the usual projection map  $\pi  : \overline{\tH} \rightarrow X_{\Gamma}(\C)$. Recall, the map $\pi$ is unramified outside the elliptic points and the 
set of cusps $\partial(X_{\Gamma})$. Both these sets are finite.

\subsection{The rational structure of the curve $X_0(N)$ defined over rational}
There is a smooth projective curve $X_0(N)$ defined over 
$\Q$ for which the space $\Gamma_0(N ) \backslash \overline{\tH}$ is canonically identified 
with the set of $\C$-points of the projective curve $X_0(N)$.
We are interested to understand the $\Q$-structure of the compactified modular curve $X_0(N)$.
\subsection{Classical modular symbols}
\label{Modularsymbols}
Recall the following fundamental theorem of Manin~\cite{MR0314846}. 
\begin{thm}
\label{MaintheoremManin}
For $\alpha \in \overline{\tH}$, consider the map $c: \Gamma \rightarrow \HH_1(X_0(N), \Z)$
defined by $$c(g)=\{\alpha, g \alpha\}.$$
The map $c$ is a surjective group homomorphism which does not depend on the choice of point $\alpha$.
 The kernel of this homomorphism is generated by 
 \begin{enumerate}
\item
 the commutator, 
 \item
  the elliptic elements, 
  \item
   the parabolic elements
\end{enumerate}
of the congruence subgroup $\Gamma$. 
\end{thm}
In particular, the above theorem implies that $\{\alpha, g\alpha\}=0$ for all $\alpha \in \sP^1(\Q)$ 
and $g \in \Gamma$. 

\subsection{The Manin map}
  Let $T, S$ be the matrices $\left(\begin{smallmatrix}
1 & 1\\
0 & 1\\
\end{smallmatrix}\right),\left(\begin{smallmatrix}
0 & -1\\
1 & 0\\
\end{smallmatrix}\right)$
and  $R=ST$ be the matrix $\left(\begin{smallmatrix}
0 & -1\\
1 & 1\\
\end{smallmatrix}\right)$. 
The modular group $\SL_2(\Z)$ is generated by 
$S$ and $T$. 
\begin{thm}[Manin]\cite{MR0314846}
Let 
$$\xi : \mathrm{SL}_2(\Z) \rightarrow \HH_1(X_0(pq),\partial(X_0(pq)),\Z)$$
be the map that takes a matrix $g \in \mathrm{SL}_2(\Z)$ to 
the class in $\HH_1(X_0(pq),\partial(X_0(pq)), \Z)$ of the image in 
$X_0(pq)$ of the geodesic in $\tH \cup \sP^1(\Q)$ joining $g.0$ and $g.\infty$. 

\begin{itemize}
 \item 
The map $\xi$ is surjective. 
\item
For all $g \in \Gamma_0(pq) \backslash \SL_2(\Z)$, 
$\xi(g)+\xi(g S)=0$ and $\xi(g)+\xi(g R)+\xi(g R^2)=0$.
\end{itemize}
\end{thm}

We have a short exact sequence,
\[
0  \rightarrow  \HH_1(X_0(pq),\Z) \rightarrow \HH_1(X_0(pq), \partial(X_0(pq)),\Z)  \rightarrow {\Z}^{\partial(X_0(pq))}\xrightarrow{\delta'} \Z \rightarrow 0.
\]
The first map is a canonical injection. The boundary map $\delta'$ takes a geodesic,
joining the cusps $r$ and $s$ to 
the formal symbol $[r]- [s]$ and the third map is the sum of the coefficients.  

\subsection{Relative homology group $ \HH_1(X_0(pq) - R \cup I, \partial(X_0(pq)), \Z)$}
\label{relative homology}
 Consider the points $i =\sqrt{-1}$
and $\rho =\frac{1+\sqrt{-3}}{2}$ on the complex upper half plane with $\nu$ the geodesic joining $i$ and $\rho$.  These are the elliptic points on the Riemann surface $X_0(pq)$. The projection map $\pi$ is 
unramified  outside cusps and elliptic points.

Say $R = \pi(\SL_2(\Z)\rho)$ and
$I = \pi(\SL_2(\Z)i)$ be the image of these two sets in $X_0 (pq)$. These two sets are disjoint.
Consider now the relative homology group $\HH_1(Y_0(pq), R \cup I, \Z)$.
 For $g \in  \SL_2(\Z)$, let
$[g]_*$ be the class of $\pi(g \nu)$ in the relative homology group $\HH_1(Y_0(pq), R \cup I, \Z)$. Let $\rho^* = -\overline{\rho}$ 팫 be another
point on the boundary of the fundamental domain. The homology groups $\HH_1(Y_0(pq), \Z)$ are subgroups of $\HH_1(Y_0(pq), R \cup I, \Z)$. 
Suppose $z_0 \in  \tH$ be such that $|z_0| = 1$ and $\frac{-1}{2} < Re(z_0 ) < 1$. Let $\gamma$ be the union of the geodesic in $\tH \cup \sP^1(\Q)$
joining $0$ and $z_0$ and $z_0$ and $i \infty$. For $g \in \Gamma_0(pq) \backslash \SL_2(\Z)$, let $[g]^*$ be the class of $\pi(g \gamma)$
in $\HH_1(X_0(pq)- R \cup I, \partial(X_0(pq)), \Z)$. 

We have an intersection pairing 
\[
\circ : \HH_1(X_0(pq) - R \cup I, \partial(X_0(pq)), \Z) \times \HH_1(Y_0(pq), R \cup I, \Z) \rightarrow \Z.
\]
 Recall the following results of Merel [Prop.\,1, Cor.\,1, \cite{MR1363498}].
\begin{prop}\cite{MR1405312}
\label{rho-value}
For $g,h \in \Gamma_0(pq) \backslash \SL_2(\Z)$, we have 
\[
 [g]^* \circ [h]_*=1 
\]
if $\Gamma_0(pq) g = \Gamma_0(pq) h$ and 
 \[
 [g]^* \circ [h]_*=0
\]
otherwise.
\end{prop}
\begin{cor}
\label{rho-int}
 The homomorphism of groups $\Z^{\Gamma_0(pq)\backslash \SL_2(\Z)} \rightarrow \HH_1(Y_0 (pq), R \cup I, \Z)$ induced by the map
$$\xi_0 (\sum_g \mu_g g) = \sum_g \mu_g[g]_*$$
is an isomorphism.
\end{cor}
 The following important property [Cor.\, 3, \cite{MR1363498}] of the intersection pairing will be used later.
\begin{cor}
 For $g \in \Gamma_0(pq) \backslash \SL_2(\Z)$, let $\sum_h \mu_h h \in  {\Z}^{\Gamma_0(pq) \backslash \SL_2(\Z)}$ be such that $\sum_h \mu_h [h]_*$ is the image
of an element of $\HH_1(Y_0(pq), \Z)$ under the canonical injection. We have
\[
[g]^* \circ (\sum_h \mu_h [h]_*)=\mu_g. 
\]
\end{cor}
We have a short exact sequence,
\[
0  \rightarrow  \HH_1(X_0(pq)-R \cup I ,\Z) \rightarrow \HH_1(X_0(pq)-R \cup I , \partial(X_0(pq)),\Z)  \rightarrow {\Z}^{\{\partial(X_0(pq))\}} \xrightarrow{\delta}  \Z \rightarrow 0.
\]
The boundary map $\delta$ takes a geodesic,
joining the cusps $r$ and $s$ to the formal symbol $[r]- [s]$. Note that
$\delta'(\xi(g))=\delta([g]^*)$ for all $g \in \SL_2(\Z)$.

Recall,  we have a canonical bijection $\Gamma_0(pq)\backslash \mathrm{SL}_2(\Z) \cong
{\sP}^1(\Z/pq\Z)$ given 
by $\left(\begin{array}{cc}
a & b\\
c & d\\
\end{array}\right) \rightarrow (c,d)$. Say $\alpha_k , \beta_r$ and $\gamma_s$ are the matrices 
 $\left(\begin{smallmatrix}
0 & -1\\
1 & k\\
\end{smallmatrix}\right)$,  $\left(\begin{smallmatrix}
-1 & -r\\
p & rp-1\\
\end{smallmatrix}\right)$ and $\left(\begin{smallmatrix}
-1 & -s\\
q & sq-1\\
\end{smallmatrix}\right)$ respectively.
We explicitly write down the elements of $\sP^1(\Z/pq\Z)$ as a set
\[
 \{(1, k), (1, tp), (1, t'q), (p, q), (q, p), (tp, 1), (t' q, 1), (1, 0), (0, 1)\}
\]
with $k \in  (\Z/pq\Z)^* , t \in  (\Z/q\Z)^* , t' \in (\Z/p\Z)^*$ . Observe that $(p, q) = (t p, q) = (p, t'q)$ for all $t$ and $t'$ co
prime to $pq$.

\begin{lemma}
\label{explicit}
The set  $\varOmega=\{I,\alpha_{k},\beta_r, \gamma_s| 0 \leq k \leq pq-1,  0 \leq r  \leq (p-1), 0\leq s \leq (q-1) \}$
forms a complete set of  coset representatives of $\Gamma_0(pq)\backslash \mathrm{SL}_2(\Z)$. 
\end{lemma}
\begin{proof}
The  orbits $\Gamma_0(pq) \alpha_k$, $\Gamma_0(pq) \beta_{l}$ and $\Gamma_0(pq) \gamma_{m}$ are disjoint
since $ab^{-1}$ do not belong to $\Gamma_0(pq)$ for two distinct matrices $a$, $b$ from the set $\varOmega$.
There are $1+ pq + p + q =|{\sP}^1(\Z/pq\Z)|$ coset representatives.  
\end{proof}

We list different rational numbers of the form $\alpha (0)$ and  $\alpha( \infty)$ with 
$\alpha \in \varOmega$ as equivalence classes of cusps as follows:
\[
\begin{tabular}{|l|l|l|}
\hline
$0$ & $\frac{1}{p}$ & $\frac{1}{q}$ \\ \hline
$\frac{-l}{lp-1}$, $(lp-1,q)=1$ &   $\frac{-1}{k}$, $(k,p) >1$ & $\frac{-1}{k}$, $(k,q) >1$\\ \hline
$\frac{-m}{mq-1}$, $(mq-1,p)=1$& $\frac{-m}{mq-1}$, $(mq-1,p) >1$  & $\frac{-l}{lp-1}$, $(lp-1,q) >1$ \\ \hline
\end{tabular}.
\]
\subsection{Manin-Drinfeld theorem}
\label{ManinDrinfeld}
Following ~\cite{MR1363488}, we briefly recall the statement of the Manin-Drinfeld 
theorem. 
\begin{thm}[ Manin-Drinfeld]
\label{ManinDrinfeld}
 \cite{MR0318157}
For a congruence subgroup $\Gamma$ and any two cusps $\alpha$, $\beta$ in $\sP^1(\Q)$, 
the path 
\[
\{ \alpha, \beta \} \in \HH_1(X_{\Gamma}, \Q).
\]
\end{thm}
This theorem can be reformulated in terms of divisor classes on the Riemann surface. 
\begin{thm}

Let $a=\sum_i m_i P_i$ be a divisor of degree zero on $X$. Then $a$ is a divisor of a rational 
function if and only if there exist a cycle $\sigma \in \HH_1(X_{\Gamma}, \Z)$ such that 
\[
\int_a \omega=\sum_i m_i \int_{P_0}^{P_i} \omega=\int_{\sigma} \omega
\] 
for every $\omega \in \HH^0(X_{\Gamma}, \Omega_{X_{\Gamma}})$. 
\end{thm}
As a corollary, we notice that $\{x, y\} \in \HH_1(X_{\Gamma}, \Q)$ if and if there is a positive integer 
$m$ such that $m(\pi_{\Gamma}(x)-\pi_{\Gamma}(y))$ is a divisor of a function. In other words, 
the degree zero divisors supported on the cusps are of finite order in the divisor class 
group. Manin-Drinfeld proved it using the extended action of the usual Hecke operators.
In particular, it says that $\{0, \infty\} \in \HH_1(X_{\Gamma}, \Q)$ although $0$ and $\infty$
are two inequivalent cusps of $X_{\Gamma}$. 
In \cite{MR0364259}, Ogg constructed certain modular function $X_0(pq)$ whose divisors 
coincide with degree zero divisors on the modular curves. 
\section{Eisenstein series for \texorpdfstring{$\Gamma_0(pq)$}{X} with integer coefficients}
\label{Eisensteinseries}
Let $\sigma_1(n)$ denote 
the sum of the positive divisors of $n$.
We consider the series $$E'_2(z)=1-24(\sum_n \sigma_1(n) e^{2 \pi i n z}).$$
Let $\Delta$ be the Ramanujan's cusp form of weight $12$.
For all $N \in \N$, the function 
$z \rightarrow \frac{\Delta(Nz)}{\Delta(z)}$ is a function on $\tH$ invariant under $\Gamma_0(N)$.  
The logarithmic differential of this function is 
$2 \pi i E_N(z) dz$ and  $E_N$ is a classical holomorphic modular form of weight  two for 
$\Gamma_0(N)$ with constant 
term $N-1$. The differential form $E_N(z) dz$ is a differential form of third kind on $X_0(N)$. The 
periods [\S\,\ref{Peiods}] of these differential forms are in $\Z$. 

By  [\citep{MR2112196}, Thm.\, 4.6.2], the set $\tE_{pq}=\{E_p, E_q, E_{pq} \}$ 
is a basis of $E_2(\Gamma_0(pq))$. 

\begin{lemma}
\label{cuspspq}
The cusps $\partial(X_0(pq))$ can be identified with the set $\{ 0, \infty, 
\frac{1}{p}, \frac{1}{q} \}$. 
\end{lemma}
\begin{proof}
If $\frac{a}{c}$ and $\frac{a^{\prime}}{c^{\prime}}$ are in $\sP^1(\Q)$, then 
 $\Gamma_0(pq) \frac{a}{c} = \Gamma_0(pq) \frac{a^{\prime}}{c^{\prime}}$  
$\Longleftrightarrow$  $\left(\begin{array}{c} a y \\ c \\ \end{array}\right)  \equiv \left(\begin{array}{c} a^{\prime} +
jc^{\prime} \\ c^{\prime}y\\ \end{array}\right)$  $\pmod{pq}$, for some $j$ and $y$ such that $\mathrm{gcd}(y,pq)=1$
[cf.~\cite{MR2112196}, p. 99]. A small check shows that the orbits 
$\Gamma_0(pq)0$, $\Gamma_0(pq)\infty$, $\Gamma_0(pq)\frac{1}{p}$ and $\Gamma_0(pq)\frac{1}{q}$ 
are disjoint.
\end{proof}
Let $\mathrm{Div}^{0}(X_0(pq), \partial(X_0(pq)),\Z)$ be the group of degree zero divisors supported on cusps.
For all cusps $x$, let $e_{\Gamma_0(pq)}(x)$ denote the ramification index of $x$ over $\SL_2(\Z) \backslash \tH \cup \sP^1(\Q)$ and 
\[
 r_{\Gamma_0(pq)}(x)=e_{\Gamma_0(pq)} (x)a_0(E[x]).
\]
By [\cite{MR670070}, p.\, 23], there is a canonical isomorphism $\delta : E_2(\Gamma_0(pq)) \rightarrow \mathrm{Div}^0(X_0 (pq), \partial(X_0(pq)), \Z)$ that takes the Eisenstein series $E$ to the divisor
\begin{equation}
\label{divisor}
\delta(E)=\sum_{x \in \Gamma_0(pq) \backslash \sP^1(\Q)} r_{\Gamma_0(pq)} (x)[x].
\end{equation}
Hence,  the Eisenstein element is related to the Eisenstein series by 
the boundary map. In Prop.~\ref{boundact}, we prove that the boundary of Eisenstein element is indeed the boundary of Eisenstein series. 
By [\cite{MR800251}, p. 538], we see that  
\begin{equation*}
 e_{\Gamma_0(pq)}(x) =
\begin{cases}
q & \text{if $x=\frac{1}{p}$} \\
p & \text{if $x=\frac{1}{q}$}\\
1 & \text{if $x= \infty$}\\
pq & \text{if $x = 0.$} 
\end{cases} 
\end{equation*}
Since $\displaystyle \sum_{x \in \partial(X_0(pq))}e_{\Gamma_0(pq)}(x)a_0(E[x])=0$, we write the corresponding degree zero divisor as 
$$\delta(E)=  a_0(E) (\{ \infty \} - \{0 \}) + q a_0(E[\frac{1}{p}]) 
(\{ \frac{1}{p}\} - \{ 0 \}) + p a_0(E[\frac{1}{q}]) (\{ \frac{1}{q} \} - \{ 0 \}).$$
\subsection{Period Homomorphisms}
\label{Peiods}
We now define the period homomorphisms for the differential forms of third kind. 
\begin{definition}[Period homomorphism]
For $E_N \in \tE_{pq}$,  the differential forms $E_N(z)dz$ are of third kind on the Riemann surface
$X_0(pq)$ but of first kind on the non-compact Riemann surface $Y_0(N)$.
 For any $z_0 \in \tH$ and $\gamma \in \Gamma_0(pq)$, let $c(\gamma)$ 
be the class in $\HH_1(Y_0(pq),\Z)$ of the image in $Y_0(pq)$ of the geodesic in $\tH$ joining 
$z_0$ and $\gamma(z_0)$. 
That  the class is non-zero follows from Thm.~\ref{MaintheoremManin}.
This class is independent of the choice of $z_0 \in \tH$ and let $\pi_{E_N}(\gamma)=\int_{c(\gamma)} E_N(z)dz$. 
The map $\pi_{E_N}:\Gamma_0(pq) \rightarrow \Z$
is the ``period'' homomorphism of $E_N$.
\end{definition}
Let $\overline{B}_1(x)$ be the first Bernoulli's polynomial of period $1$ defined by 
\[
\overline{B}_1(0)=0, \overline{B}_1(x)=x-\frac{1}{2}
\]
if $x \in (0, 1)$. For any two integers $u$ and $v$  with $v \geq 1$, we define the Dedekind sum 
by the formula:
\[
S(u,v)=\sum_{t=1}^{v-1} \overline{B}_1(\frac{tu}{v}) \overline{B}_1(\frac{u}{v}). 
\]
Recall some well-known properties of the period mapping $\pi_{E_N}$ [cf.~\cite{MR553997}, p.\,10, \cite{MR1405312}, p.\,14] for the Eisenstein series $E_N \in  \tE_{pq}$.
\begin{prop}
\label{value-pi}
Let $\gamma=\left(\begin{smallmatrix}
a & b\\
c & d\\
\end{smallmatrix}\right)$ be an element of $\Gamma_0(pq)$. 
\begin{enumerate}
 \item 
$\pi_{E_N}$ is a homomorphism $\Gamma_0(pq) \rightarrow \Z$. 
\item 
Consider the number $\mu=gcd(N-1,12)$, the image of $\pi_{E_N}$ lies in $\mu \Z$.
\item 
\begin{equation*}
\pi_{E_N}(\gamma) = 
\begin{cases}
\frac{a+d}{c} (N-1)+12 sgn(c)(S(d,|c|)-S(d,\frac{|c|}{N})) & \text{if $c \neq 0$,} \\
\frac{b}{d} (N-1) & \text{if $c=0$.}
\end{cases}
\end{equation*}
\item 
\[
\pi_{E_N}(\gamma)=\pi_{E_N}(\left(\begin{smallmatrix}
d & \frac{c}{N}\\
Nb & a\\
\end{smallmatrix}\right))
\]
\end{enumerate}
\end{prop}

\section{Eisenstein elements}
\label{Eisensteinelements}
Following ~\cite{MR1405312} and ~\cite{MR1240651}, recall the concept of Eisenstein 
elements of the space of modular symbols. 
For any natural number $M>4$, the congruence subgroup $\Gamma_0(M)$ is the subgroup of $\mathrm{SL}_2(\Z)$ consisting of all matrices 
$\left(\begin{smallmatrix}
a & b\\
c & d\\
\end{smallmatrix}\right)$ such that $M \mid c$.
The congruence subgroup $\Gamma_0(M)$
acts on the upper half plane $\tH$ in
the usual way. 
The quotient space $\Gamma_0(M) \backslash \tH$
is denoted by $Y_0(M)$.  
Apriori,  these are all  Riemann surfaces and hence  algebraic curves defined over $\C$. 
There are  models of these algebraic curve defined over $\Q$ and they parametrize elliptic curves with 
cyclic subgroups of order $M$.  Let $X_0(M)$ be the compactification of the Riemann surface $Y_0(M)$ obtained by adjoining the set 
of cusps $\partial(X_0(M))=\Gamma_0(M) \backslash \sP^1(\Q)$. 
\begin{definition}
[Eisenstein elements]
Let $\pi_{E_N}:\HH_1(Y_0(pq),\Z) \rightarrow \Z$
be the ``period'' homomorphism of $E_N$ [\S~\ref{Peiods}].
The intersection pairing $\circ$ \cite{MR1240651} induces a perfect, bilinear pairing 
\[
 \HH_1(X_0(pq),\partial(X_0(pq),\Z) \times \HH_1(Y_0(pq),\Z) \rightarrow \Z. 
\]
Since $\circ$ is a non-degenerate bilinear pairing, 
there is an unique element $\sE_{E_N} \in \HH_1(X_0(pq),\partial(X_0(pq)),\Z)$ such that $\sE_{E_N} \circ c=\pi_{E_N}(c).$ 
The modular symbol $\sE_{E_N}$ is the {\it Eisenstein element} corresponding to the Eisenstein series 
$E_N$.
\end{definition}
We intersect with the congruence 
subgroup $\Gamma(2)$ to ensure that the Manin maps become bijective (rather than only 
surjective), compute the Eisenstein elements  for these modular curves, calculate the boundary 
and show that the these boundaries  coincide with the original Eisenstein elements.  
In the case of $\Gamma_0(p^2)$, although it is difficult to find the Fourier 
expansion of modular forms at different cusp but fortunately for all $g \in \Gamma_0(p)$
the matrices 
$g \left(\begin{smallmatrix}
1 & 1\\
0 & 1\\
\end{smallmatrix}\right) g^{-1}$ belong to $\Gamma_0(p^2)$ and hence it was easier to tackle the explicit coset representatives. Unfortunately, 
for $N=pq$ or $N=p^3$ these are no longer true.

To get around this problem for the congruence subgroup $\Gamma_0(pq)$ with $p$ and $q$ distinct 
primes, we use the relative homology group $\HH_1(X_0(pq), R \cup I, \Z)$. For these relative 
homology groups, the associated Manin maps are bijective  and the push forward of these Eisenstein elements 
inside the original modular curve turn out to have same boundary as the original Eisenstein 
elements. We consider three different homology groups in these paper and in particular the study
of the relative homology group $\HH_1(X_0(N), R \cup I, \Z)$  to determine the Eisenstein element is a new idea that we wish to propose 
in this article. That these 
relative homology groups should be useful in the study of modular symbol are discovered 
by Merel.

\begin{definition}[Almost eisenstein elements]
For $N \in \{p,q,pq\}$, the differential form $E_N(z)dz$ is of first kind on the Riemann surface $Y_0(pq) $. Since $\circ$ is a
non-degenerate bilinear pairing, there is an unique element $\sE'_{E_N} \in  \HH_1(X_0(pq) - R \cup I, \partial(X_0(pq)), \Z)$ such that
$\sE'_{E_N} \circ c = \pi_{E_N}(c)$ for all $c \in  \HH_1(Y_0 (pq), R \cup I, \Z)$. We call $\sE'_{E_N}$ the {\it almost Eisenstein element} corresponding
to the Eisenstein series $E_N$.
\end{definition}

\section{Even Eisenstein elements}
\label{Even}
\subsection{Simply connected Riemann surface of genus zero with three marked points}
Recall, there is only one simply connected (genus zero) compact Riemann surface up to conformal bijections: namely  the Riemann sphere  or  the projective complex plane $\sP^1(\C)$. 
A theorem of Belyi states that every (compact, connected, non- singular) algebraic curve $X$ 
has a model defined over $\overline{\Q}$  if and only if it admits a map to $\sP^1(\C)$ branched over three points.

Consider the subgroup $\Gamma(2)$ of $\mathrm{SL}_2(\Z)$ consisting of all matrices 
which are identity modulo the reduction map modulo $2$.  The Riemann 
surface $\Gamma(2) \mod \overline{\tH}$ is a Riemann surface of genus zero, denoted by $X(2)$. Hence, it can be 
identified with $\sP^1(\C)$.

The subgroup $\Gamma(2)$ has three cusps $\Gamma(2) 0$,  $\Gamma(2) 1$ and 
 $\Gamma(2) \infty$. Hence, $\Gamma(2) \backslash \overline{\tH}$ become the simply 
 connected Riemann surface $\sP^1(\C)$ with the three marked points 
 $\Gamma(2) 0$,  $\Gamma(2) 1$ and 
 $\Gamma(2) \infty$ given by respective cusps. The modular curve $X_0(pq)$ has no 
 obvious morphism to $X(2)$. Hence, we consider the modular curve $X_{\Gamma}$ [\S
 \ref{bijectiveManin}]. There are two obvious maps $\pi, \pi'$ from $X_{\Gamma}$ to the 
 compact Riemann surface $X_0(pq)$. 

\subsection{Modular curves with bijective Manin maps}
\label{bijectiveManin}
For the congruence subgroup  $\Gamma=\Gamma_0(pq)
\cap \Gamma(2)$, consider the  compactified modular curve $X_\Gamma=\Gamma \backslash \tH \cup
\sP^1(\Q)$ and let 
$ \pi_{\Gamma} :\tH \cup \sP^1(\Q) \rightarrow X_\Gamma$ be the canonical surjection. 

Let $\pi_0:\Gamma \backslash  \tH \cup \sP^1(\Q) \rightarrow \Gamma(2) \backslash \tH \cup \sP^1(\Q) $ be the map $\pi_0(\Gamma z)=\Gamma(2)z$.
The compact Riemann surface  $X(2)$ contains three cusps $\Gamma(2)1$, $\Gamma(2)0$, $\Gamma(2)\infty$.
Let $P_{-}=\pi_0^{-1}(\Gamma(2)1)$ and $P_{+}$ be the union 
of two sets $\pi_0^{-1}(\Gamma(2)0)$ and $\pi_0^{-1}(\Gamma(2)\infty)$.
Consider now the Riemann surface $X_{\Gamma}$ with boundary $P_{+}$ and $P_{-}$.

Let $\delta_r$ be $1$ or $0$ depending on $r$ is odd or even. 
For any integer $k$, let $s_k=(k+(\delta_k-1)pq)$ be an odd integer. Say $l$ and $m$ be two unique integers such
that $lq + mp \equiv  1 \pmod{pq}$ with $1 \leq l \leq (p - 1)$ and $1 \leq m \leq (q - 1)$. The matrices
$\alpha'_{pq}=\left(\begin{smallmatrix}
pq  &  pq-1\\ 
pq+1 &  pq\\
\end{smallmatrix}\right)$,  $\alpha'_k=\left(\begin{smallmatrix}
s_k(pq)^2   &  s_kpq-1\\ 
s_kpq+1 &  s_k\\
\end{smallmatrix}\right)$ $\beta^{\prime}_r=\left(\begin{smallmatrix}
-1 & -(r+\delta_r q)\\
p+pq & -1 +(r+\delta_r q)(p+pq)\\
\end{smallmatrix}\right)$ and\\$\gamma'_s=\left(\begin{smallmatrix}
-1 & -(s+\delta_s pq) \\
q+pq & -1 +(s+\delta_s pq)(q+pq)\\
\end{smallmatrix}\right)$ are useful to calculate the boundaries of the Eisenstein elements.
\begin{lemma}
\label{coset}
The set
$\Delta=\{ I, \alpha^{\prime}_{k}, \beta^{\prime}_r, \gamma^{\prime}_s  | 0 \leq k \leq (pq-1), 0 \leq r \leq (q-1), 0 \leq s \leq (p-1)\}$  $\subset$ $\Gamma(2)$ 
forms an explicit set of coset representatives of ${\sP}^1(\Z/pq\Z).$
\end{lemma}
\begin{proof}
An easy check shows that the orbits 
$\Gamma_0(pq) \alpha'_{k}$, $\Gamma_0(pq) \beta'_r$ and $\Gamma_0(pq)\gamma'_s$ 
are disjoint. Since $|\sP^1(\Z/pq\Z)|=pq+p+q+1$, the result follows.
\end{proof}

The coset representatives in the above lemma are chosen such that $\Gamma_0(pq) \beta_r$ = $\Gamma_0(pq) \beta^{\prime}_r$ and 
$\Gamma_0(pq) \gamma_s$=$\Gamma_0(pq) \gamma^{\prime}_s.$

\begin{lemma}
 $\Gamma \backslash \Gamma(2)$ is isomorphic to ${\sP}^1(\Z/pq\Z)$
\end{lemma}
\begin{proof}
 The explicit closet representatives of Lemma~\ref{coset} produce the canonical bijection. 
\end{proof}

We study the relative homology groups $\HH_1(X_\Gamma-P_{-},P_{+},\Z)$ and 
$\HH_1(X_\Gamma-P_{+},P_{-},\Z)$. The intersection pairing is a non-degenerate bilinear pairing 
$\circ:\HH_1(X_\Gamma-P_{+},P_{-},\Z) \times \HH_1(X_\Gamma-P_{-},P_{+},\Z) \rightarrow \Z.$
Recall the following two fundamental theorems from \citep{MR1405312}. For $g \in \Gamma \backslash \Gamma(2)$, let $[g]^0$ 
(respectively  $[g]_0$) be the image in $X_\Gamma$ of 
the geodesic in $\tH \cup \sP^1(\Q)$ joining $g0$ and $g\infty$ (respectively  $g1$ and $g(-1)$). 
Recall the following two fundamental theorems of~\cite{MR1405312}.
\begin{thm}[~\cite{MR1405312}]
\label{int-iso}
Let
\[
\xi_0: {\Z}^{\Gamma \backslash \Gamma(2)} \rightarrow \HH_1(X_\Gamma-P_{+},P_{-},\Z)
\]
be the map which takes $g \in \Gamma \backslash \Gamma(2)$ to the element $[g]_0$ and 
\[
\xi^0:  {\Z}^{\Gamma \backslash \Gamma(2)} \rightarrow \HH_1(X_\Gamma-P_{-},P_{+},\Z)
\]
be the map which takes $g \in \Gamma\backslash \Gamma(2)$ to the element $[g]^0$. 
The homomorphisms $\xi_0$ and $\xi^0$ are isomorphisms.
\end{thm}
\begin{theorem}[~\cite{MR1405312}]
\label{int-value}
For $g, g' \in \Gamma(2)$, we have 
\[
 [g]_0 \circ [g']^0=1 
\]
if $\Gamma g = \Gamma g'$ and 
 \[
 [g]_0 \circ [g']^0=0
\]
otherwise.
\end{theorem}
 The following two lemmas about the set $P_{-}$ are true for the congruence 
subgroup $\Gamma_0(N)$ with $N$ odd.  
\begin{lemma}
We can explicitly write the set $P_{-}$ is of the form $\Gamma \frac{x}{y}$ with $x$ and $y$ both odd. 
\end{lemma}
\begin{proof}
If possible, some element of $P_{-}$ is of the form $\Gamma \frac{x}{y}$ with $x$ and $y$ 
co-prime and $y$ even.  Consider the corresponding element in the marked simply 
connected Riemman surface $X(2)$.  The cusp $\Gamma(2) \frac{x}{y}$ is an element 
such that $y$ even and $p$ odd ($\gcd(x,y)=1$). 
First, choose $p',q'$ such that $xq'-yp'=1$ and hence $\left(\begin{smallmatrix}
x & p'\\
y & q'\\
\end{smallmatrix}\right) \in \SL_2(\Z)$. Clearly, $q'$ is odd since $y$ is even. If 
$p'$ is odd then replace the matrix   $\left(\begin{smallmatrix}
x & p'\\
y & q'\\
\end{smallmatrix}\right)$ with $\left(\begin{smallmatrix}
x & p'\\
y & q'\\
\end{smallmatrix}\right) T^{-1}$  to produce a matrix in $\Gamma(2)$  that takes $i \infty$ 
to $\frac{x}{y}$.
This is a contradiction to the fact that $\Gamma \frac{x}{y} \in P_{-}$.

If $x$ is even then the projection of $\Gamma\frac{x}{y}$ 
produces an element of $\Gamma(2)0$.  Hence, $x$ is necessarily odd. 
\end{proof}
 The following lemma is deeply 
influenced by an  important propositions of Manin [\cite{MR0314846}, Prop.\,2.2] and  Cremona [\cite{MR1628193}, Prop.\,2.\,2.\,3].
\begin{cor}
\label{Pminus}
We can explicitly write the set $P_{-}=\{\Gamma 1, \Gamma \frac{1}{pq}, \Gamma \frac{1}{p},  \Gamma \frac{1}{q}\}$
\end{cor}
\begin{proof}
Since $P_{-}=\pi_0^{-1}(\Gamma(2)1)$, we can write every element of the set $P_{-}$ as $\Gamma \theta1$ for some 
$\theta \in \Delta$ (Lemma~\ref{coset}). Let $\delta \in \{1, p,q,pq\}$, 
then every element of $P_{-}$ can be written as $\Gamma \frac{u}{v\delta}$ with 
$\gcd(u, v\delta)=1$ and $\gcd(v\delta, \frac{pq}{\delta})=1$.
Choose  an odd integer $m$ and an even integer $l$  such that $lu-mv\delta=1$. 
A calculation using matrix multiplication shows that 
 $\left(\begin{smallmatrix}
1 & 0 \\
\delta-1 & 1\\
\end{smallmatrix}\right)\left(\begin{smallmatrix}
1+c & -c \\
c& 1-c\\
\end{smallmatrix}\right)1=\frac{1}{\delta}$ and $\left(\begin{smallmatrix}
-m & u+m\\
-l & l+v\delta\\
\end{smallmatrix}\right)1=\frac{u}{v\delta}$ and hence 
$A=\left(\begin{smallmatrix}
1 & 0 \\
\delta-1 & 1\\
\end{smallmatrix}\right)\left(\begin{smallmatrix}
1+c & -c \\
c& 1-c\\
\end{smallmatrix}\right)\left(\begin{smallmatrix}
 l+v\delta & -m-u\\
l & -m\\
\end{smallmatrix}\right)$ is a matrix such that $A(\frac{u}{v\delta})=\frac{1}{\delta}$.  The matrix 
$A$ belongs to $\Gamma$ if and only if $c v\delta \equiv l' \pmod {\frac{pq}{\delta}}$. Since 
$v \delta$ is coprime to $\frac{pq}{\delta}$, there is always such $c$.
Hence, the set $P_{-}$ consists of four elements as in the statement of the Corollary. 
\end{proof}
Let $\pi, \pi':\Gamma \backslash  \overline{\tH} \rightarrow \Gamma_0(pq) \backslash \overline{\tH}$ 
be the maps $\pi(\Gamma z)=\Gamma_0(pq)z$ and 
$\pi'(\Gamma z)=\Gamma_0(pq)\frac{z+1}{2}$ respectively. 
Consider the matrix $h=\left(\begin{smallmatrix}
1 & 1\\
0 & 2\\
\end{smallmatrix}\right)$.
The morphism $\pi'$ is well-defined since the matrix  $h \gamma h^{-1}$ 
belongs to $\Gamma_0(pq)$ for all $\gamma \in \Gamma$.
The morphisms $\pi, \pi'$ together induce a map 
\[
\kappa: \C(X_{\Gamma}) \rightarrow \C(X_0(pq))
\]
between the function fields of the Riemann surfaces $X_{\Gamma}$ by $\kappa(f(z))=\frac{f(\pi(\Gamma z))^2}{f( \pi'(\Gamma z))}$. 
Recall  the description of the coordinate chart around a cusp $\Gamma x$ \cite{MR2194815} of the Riemann surface $X_{\Gamma}$. 
\begin{definition}
\label{coordinatecharts}
For a cusp $y$ of the congruence subgroup $\Gamma$, let $\Gamma_y$ be the subgroup of 
$\Gamma$ fixing the cusp $y$. 
 Let $t \in \SL_2(\R)$ be such that $t(y)=i \infty$ and $m$ be the smallest natural number such that 
$t \Gamma_{y} t^{-1}$ is generated by  the matrix $\left(\begin{smallmatrix}
1 & m \\
0 & 1\\
\end{smallmatrix}\right)$.
For the modular curve $X_{\Gamma}$, 
the local coordinate around the point $\Gamma y$ is $z \rightarrow e^{2 \pi i \frac{t(z)}{m}}$.
 
\end{definition}
\begin{example}
Let $y=\frac{1}{\delta }$ with $\delta$ one of the prime $p$ or $q$,
then $h(y)=\frac{u}{\delta}$ with $(u, pq)=1$. Choose integers $u',\delta'$ with $\delta'$ even such that 
$u \delta'-u'\delta=1$ and hence $\rho_{h(y)}=\left(\begin{smallmatrix}
\delta'& u' \\
-\delta & u\\
\end{smallmatrix}\right)$ is a matrix such that $\rho_{h(y)}(h(y))=i \infty$.  We can choose such a $\delta' \in \Z$ since $\delta$ is odd. 

A calculation using matrix multiplication shows that $\rho_{h(y)}T^{e} {\rho_{h(y)}}^{-1}=\left(\begin{smallmatrix}
1+e\delta \delta'& e (\delta')^2\\
-e\delta^2 & 1-e\delta \delta'\\
\end{smallmatrix}\right)$. Hence, the smallest possible $e$ to ensure $t T^et^{-1} \subset \Gamma_0(pq)$ 
is $\frac{pq}{\delta}$.  
\end{example}
\begin{example}
Since $\det(\rho_{h(y)} \circ h)=2$, hence $t=\left(\begin{smallmatrix}
\frac{l}{2} & 0 \\
& 1\\
\end{smallmatrix}\right)\rho_{h(y)} \circ h \in \SL_2(\R)$ 
and $t(y)=i \infty$. 
A calculation using matrices shows that $t T^{e} t^{-1}=\left(\begin{smallmatrix}
1+\frac{e\delta \delta'}{2}& \frac{e {\delta'}^2}{4}\\
-e\delta^2 & 1-\frac{e\delta \delta'}{2}\\
\end{smallmatrix}\right)$. Hence, the smallest possible $e$ to ensure $t T^e t^{-1} \subset \Gamma$ 
is $e=2\frac{pq}{\delta}$.  
\end{example}

We use the following lemma to construct differential forms of first kind on the ambient 
Riemann surface $X_{\Gamma}-P_{+}$. 
\begin{lemma}
\label{no-zero}
 Let $f: X_0(pq) \rightarrow \C$ be a rational function. The divisors of $\kappa(f)$ 
are supported on  $P_{+}$.
\end{lemma}
\begin{proof}
Suppose $f$ is a meromorphic function on the Riemann surface $X_0(pq)$. Then it is given by $\frac{g}{h}$ with $g$ and $h$ holomorphic function on the Riemann surface $X_0(pq)$. 
Every element of $P_{-}$ is of the form $\Gamma \frac{1}{\delta}$ with $\delta \mid N$. 
By [Prop.\, 4.1, p.\,44, \cite{MR1326604}], every holomorphic map on Riemann surface locally looks like 
$z \rightarrow z^n$.

Consider the morphism   $\pi'$ and the point on the modular curve $\Gamma\frac{1}{\delta}$.  
The local coordinates around the point 
$\Gamma_0(pq)0,\Gamma_0(pq)\infty$ and $\Gamma_0(pq)\frac{1}{p}$ 
are given by $q_0(z)=e^{2\pi i \frac{1}{-pqz}}, q_\infty(z)=e^{2\pi i z}$ and $q_\frac{1}{q}(z)= e^{2\pi i \frac{z}{p(-qz+1)}}$
respectively. In the modular curve $X_\Gamma$, the local coordinates around the points of $P_{-}$ are given by
$q_1(z)=e^{2\pi i \frac{1}{2pq(-z+1)}}$, $q_{\frac{1}{pq}}(z)=e^{2\pi i \frac{z}{2(-pqz+1)}}$, $q_\frac{1}{p}(z)=e^{2\pi i \frac{z}{2q(-pz+1)}}$
and $q_\frac{1}{q}(z)=e^{2\pi i \frac{z}{2p(-qz+1)}}$. 

Now around the point $\Gamma 1$ and $\Gamma \frac{1}{pq}$, we have the equalities
$q_0 \circ \pi=q_1^2,q_0 \circ \pi'=q_1^4$
and, $q_{\frac{1}{pq}} \circ \pi=q_{\frac{1}{pq}}^2, q_ {\frac{1}{pq}}\circ \pi'=q_{\frac{1}{pq}}^4.$

Let $y=\frac{1}{\delta}$ with $\delta$ one of the prime $p$ or $q$. 
The local coordinate chart around the point $\Gamma \frac{1}{\delta}$ is $z \rightarrow e^{2 \pi i \frac{\rho_{h(x)}\circ h(z)}{4e}}$.  The map $\pi'$ takes it to $e^{2 \pi i \frac{2\rho_{h(x)}(h(z))}{e}}$. For this coordinate chart the map $\pi'$ is given by $z \rightarrow z^4$. 

We now consider the map $\pi$ and $t=\left(\begin{smallmatrix}
1& 0 \\
-\delta & 1\\
\end{smallmatrix}\right)$ is a matrix such that $t(y)=i \infty$ and $e=\frac{pq}{\delta}$. 
The local coordinate around the point $\Gamma \frac{1}{\delta}$ is $z \rightarrow e^{2 \pi i \frac{t(z)}{2e}}$ and the map $\pi$ takes it to $e^{2 \pi i \frac{t(z)}{e}}$. In this coordinate chart, the map $\pi$ is given by $z \rightarrow z^2$. Hence, the function $\frac{(f \circ \pi)^2}{f \circ \pi'}$ has no zero or pole on $P_{-}$.
\end{proof}
\begin{definition}
\label{Eisensteinsec}[Even Eisenstein elements]
For $E_N \in \tE_{pq}$, let 
$\lambda_{E_N}:X_0(pq) \rightarrow \C$ be the rational function whose logarithmic differential is $2 \pi i E_N(z)dz=2\pi i \omega_{E_N}$. 
Consider the rational function $\lambda_{E_N,2}=\frac{(\lambda_{E_N} \circ \pi)^2}{\lambda_{E_N} \circ \pi'}$ on $X_\Gamma$. By Lemma 
~\ref{no-zero}, this function has no zeros and poles in $P_{-}$. Let $\kappa^*(\omega_{E_N})$ be the logarithmic differential of the function.
Let $\varphi_{E_N}(c)=\int_c \kappa^*(\omega_{E_N})$ be the corresponding ``period'' homomorphism $\HH_1(X_\Gamma-P_{+},P_{-},\Z) \rightarrow \Z$.

By the non-degeneracy of the intersection pairing, there is a unique element $\sE_{E_N}^0 \in \HH_1(X_\Gamma-P_{-},P_{+},\Z)$ 
such that $\sE^0_{E_N} \circ c =\varphi_{E_N}(c)$ for all $c \in  \HH_1(X_\Gamma-P_{+},P_{-},\Z)$. The modular symbol $\sE_{E_N}^0$ is the {\it even Eisenstein element}
corresponding to the Eisenstein series $E_N$. 
\end{definition}

For $E_N \in \tE_{pq}$, define a function $F_{E_N} : \sP^1(\Z/pq\Z) \rightarrow \Z$ by
\[
F_{E_N}(g) =\varphi_{E_N}(\xi_0(g)) =\int_{g(1)}^{g(-1)}[2E_N(z) - E_N(\frac{z+1}{2})]dz.
\]
\begin{remark}
\label{matrixmultiplication}
It is easy to see that for any $\gamma=\left(\begin{smallmatrix}
a &  b \\
c & d\\
\end{smallmatrix}\right) \in \Gamma(2), h\gamma h^{-1}=\left(\begin{smallmatrix}
a+c &  \frac{b+d-a-c}{2} \\
2c & d-c\\
\end{smallmatrix}\right) \in \SL_2(\Z)$. 
\end{remark}

For any matrix $\gamma \in \Gamma$,  consider the rational number $P_N(\gamma)=\frac{2 \pi_{E_N}(\gamma)-\pi_{E_N}(h\gamma h^{-1})}{12}, t(\gamma)=b+d-a-c$ and $s(\gamma)=a+c$. 

\begin{lemma}
\label{divisibleinteger}
For $\gamma=\left(\begin{smallmatrix}
a & b \\
c & d\\
\end{smallmatrix}\right) \in \Gamma$  with $c \neq 0$,  
\[
P_N(\gamma)=sgn(t(\gamma))[2 (S(s(\gamma),|t(\gamma)|pq)-S(s(\gamma), |t(\gamma)|))
\]
\[
-S(s(\gamma),|\frac{t(\gamma)}{2}|pq)+S(s(\gamma), \frac{|t(\gamma)|}{2})].
\]
In particular, $P_N(\gamma) \in \Z$ for all $\gamma \in \Gamma$.
\end{lemma}
\begin{proof}
Recall the properties of  period homomorphism [cf.\,Prop.\,\ref{value-pi}]. 
We  calculate the corresponding periods 
\[
\pi_{E_N}(\gamma)=\pi_E(T \gamma T^{-1})=\pi_{E_N}(\left(\begin{smallmatrix}
a+c & -(a+c)+b+d \\
c & -c+d\\
\end{smallmatrix}\right))=\pi_{E_N}(\left(\begin{smallmatrix}
a+c & -(a+c)+b+d \\
c & -c+d\\
\end{smallmatrix}\right))=\pi_{E_N}(\left(\begin{smallmatrix}
d-c & \frac{c}{N} \\
t(\gamma) N & a+c\\
\end{smallmatrix}\right)).
\]
By [cf.\,Prop.\,\ref{value-pi}], we have 
\[
\pi_{E_N}(\gamma)=\frac{a+d}{t(\gamma)N}(N-1) +12 sgn(t(\gamma))[ (S(s(\gamma),|t(\gamma)|N)-S(s(\gamma), |t(\gamma)|]
\]
Similarly,  
\[
\pi_{E_N}(h\gamma h^{-1})=\pi_{E_N}(\left(\begin{smallmatrix}
a+c &  \frac{b+d-a-c}{2} \\
2c & d-c\\
\end{smallmatrix}\right)=\pi_{E_N}(\left(\begin{smallmatrix}
d-c &  \frac{2c}{N} \\
\frac{t(\gamma)N}{2} N& a+c\\
\end{smallmatrix}\right)
\]
\[
=\frac{2(a+d)}{t(\gamma)N}(N-1) +12 sgn(t(\gamma)) [S(s(\gamma),\frac{|t(\gamma)|}{2}N)-S(s(\gamma), \frac{|t(\gamma)|}{2}))]
\]
Hence, we deduce the formula as in the statement. From the formula, it is easy to see that 
 $P_N(\gamma) \in \Z$ for all $\gamma \in \Gamma$. 
\end{proof}

Let $x$ be one of the prime $p$ or $q$. 
Choose integers 
$s, s'$ and $l, l'$ such that $l(s_kx+2)-2spq=1$ and $l's_kx-2s'\frac{pq}{x}=1$. 
 Let $\gamma_1^{x,k}=\left(\begin{smallmatrix}
1+4spq &  -2l \\
- 4s(s_kx+2)pq& 1+4spq\\
\end{smallmatrix}\right)$ 
and $\gamma^{x,k}_2=\left(\begin{smallmatrix}
1+4s'\frac{pq}{x} &  -2l' \\
- 4s'(s_k)pq& 1+4s'\frac{pq}{x}\\
\end{smallmatrix}\right)$ be two matrices in  $\Gamma$. 
Since the integers $l$ and $l'$ are necessarily odd, we have 
 $\gamma_1^{x,k}(\frac{1}{s_kx+2})=-\frac{1}{s_kx+2}$ and $\gamma^{x,k}_2(\frac{1}{s_kx})=-\frac{1}{s_kx}$.  
 
 Using the formula of Lemma~\ref{divisibleinteger}, we deduce that 
 \[
 s(\gamma_1^{x,k})=1-4spq(1+s_kx), t(\gamma_1^{x,k})=-2(l-2s(s_kx+2)pq)
 \]
  and 
 \[
 s(\gamma_2^{x,k})=1-4s'pq(s_k-\frac{1}{x}), t(\gamma_2^{x,k})=-2(l'-2s's_kpq).
 \]
 We can now calculate $P_N(\gamma_1^{x,k}), P_N(\gamma^{x,k}_2)$  
 using Lemma~\ref{divisibleinteger}. 

\begin{prop}
\label{exceptional}
\begin{equation*}
F_{E_N}(g) =
\begin{cases}
12(S(r,N)-2S(r,2N)) & \text{if $g=(r-1,r+1)$,} \\
6(P_N(\gamma_1^{x,k})-P_N(\gamma^{x,k}_2)) &  \text{if $g=(1+kx,1)$ or $g=(-1-kx,1)$,} \\
-6(P_N(\gamma_1^{x,k})-P_N(\gamma^{x,k}_2)) &  \text{if $g=(1,-1-kx)$ or $g=(1,1+kx)$,} \\
0 &  \text{if $g=(\pm 1,1)$.} \\
\end{cases}
\end{equation*}
\end{prop}

\begin{proof}
If $g = (r - 1, r + 1)$ and $E_N \in \tE_{pq}$, we get [\cite{MR1405312}, p.\,18]
\[
F_{E_N}(g) = \varphi_{E_N}(\xi_0(g)) = 12(S(r, N ) - 2S(r, 2N )).
\]
We proceed to find the value of the integrals in the remaining cases. 
The differential form $k^*(\omega_{E_N})$ is of first kind on the Riemann surface $X_{\Gamma}-P_{+}$.
We also note that $g =  (\pm 1 ,1),  ( \pm1  \pm kx, 1)$ or  $(1, \pm 1 \pm kx) $ 
with $x$ one of the prime $p$ or $q$, then we can't write it as $(r-1,r+1)$.

Since all the Fourier coefficients of the Eisenstein series are real valued, so an argument similar to [\cite{MR1405312}, p.\,19] shows that 
$F_{E_N}(s_kx+1,1)=F_{E_N}(-s_kx-1,1)$. Consider the path 
$
\{\frac{1}{s_kx+2},-\frac{1}{s_kx+2}\}=\{\frac{1}{s_kx+2},\frac{1}{s_kx}\}
+\{\frac{1}{s_kx},\frac{-1}{s_kx}\}
+\{\frac{-1}{s_kx},\frac{-1}{s_kx+2}\}.
$
The rational number $\frac{1}{s_kx}$ correspond to a point of $P_{-}$ in the Riemann surface $X_{\Gamma}$. 
The differential form $k^*{\omega_{E_N}}$ has no zeros and poles on $P_{-}$. 
We deduce that 
$$
\int_{\frac{1}{s_kx+2}}^{-\frac{1}{s_kx+2}} k^*(\omega_{E_N})= \int_{\frac{1}{s_kx+2}}^{\frac{1}{s_kx}} k^*(\omega_{E_N})
+\int_{\frac{1}{s_kx}}^{\frac{-1}{s_kx}} k^*(\omega_{E_N})
+\int_{\frac{-1}{s_kx}}^{\frac{-1}{s_kx+2}} k^*(\omega_{E_N})=2 F_N(s_kx+1,1)+\int_{\frac{1}{s_kx}}^{\frac{-1}{s_kx}} k^*(\omega_{E_N}). 
$$
Let $\gamma_1^{x,k}$ and $\gamma^{x,k}_2$ be two matrices in  $\Gamma$ such that  $\gamma_1^{x,k}(\frac{1}{s_kx+2})=-\frac{1}{s_kx+2}$ and $\gamma^{x,k}_2(\frac{1}{s_kx})=-\frac{1}{s_kx}$. 
We deduce that $2 F_N(s_kx+1,1)=\int_{\frac{1}{s_kx+2}}^{\gamma_1^{x,k}(\frac{1}{s_kx+2})} k^*(\omega_{E_N})
-\int_{\frac{1}{s_kx}}^{\gamma^{x,k}_2(\frac{1}{s_kx})} k^*(\omega_{E_N})$. 

We now prove  that the $\int_{\frac{1}{s_kx}}^{\gamma^{x,k}_2(\frac{1}{s_kx})} k^*(\omega_{E_N})$ 
is independent of the choice of the matrices $\gamma^{x,k}_2 \in \Gamma$
that take $\frac{1}{s_kx}$ to $-\frac{1}{s_kx}$.  If possible, $\gamma^{x,k}_2$ and 
${\gamma'}^{x,k}_2$ be two matrices such that  $\gamma^{x,k}_2(\frac{1}{s_kx})=
{\gamma'}^{x,k}_2(\frac{1}{s_kx})=-\frac{1}{s_kx}$.  Since  $\gamma^{x,k}_2 \in \Gamma$, the integral $\varphi_{E_N}(\gamma^{x,k}_2)=\int_{\frac{1}{s_kx}}^{\gamma^{x,k}_2(\frac{1}{s_kx})} k^*(\omega_{E_N})$ is independent of the choice of any point in $\tH \cup \{-1\}$,  hence by 
replacing  $\frac{1}{s_kx}$ with $(\gamma^{x,k}_2)^{-1} ({\gamma'}^{x,k}_2)\frac{1}{s_kx}$, 
we get the above integral is same as $\int_{\frac{1}{s_kx}}^{{\gamma'}^{x,k}_2(\frac{1}{s_kx})} k^*(\omega_{E_N})$ and the integral is independent of the choice of exceptional matrices. 
Similarly, we can prove that $\int_{\frac{1}{s_kx+2}}^{\gamma^{x,k}(\frac{1}{s_kx+2})}k^*(\omega_{E_N})$ 
is also independent of the choice of the matrices that take $\frac{1}{s_kx+2}$ 
to $-\frac{1}{s_kx+2}$.  Since we have already written down two matrices 
$\gamma_1^{x,k}$ and $\gamma^{x,k}_2$ in  $\Gamma$ such that  $\gamma_1^{x,k}(\frac{1}{s_kx+2})=-\frac{1}{s_kx+2}$ and $\gamma^{x,k}_2(\frac{1}{s_kx})=-\frac{1}{s_kx}$, we use 
these matrices to find those integrals.

The above calculation shows that
\[
2 \pi_{E_N}(\gamma_1^{x,k})-\pi_{E_N}(h \gamma_1^{x,k} h^{-1})=2 F_N(s_kx+1,1)+2 \pi_{E_N}(\gamma^{x,k}_2)
-\pi_{E_N}(h \gamma^{x,k}_2h^{-1}).
\]
Hence, we get $$F_{E_N}(s_kx+1,1)=\frac{2 \pi_{E_N}(\gamma_1^{x,k})-\pi_{E_N}(h \gamma_1^{x,k} h^{-1})-2 \pi_E(\gamma^{x,k}_2)+\pi_E(h \gamma^{x,k}_2 h^{-1})}{2}=6(P_N(\gamma^{x,k})-P_N(\gamma^{x,k}_2)).$$
Since  $F_{E_N}((1+s_kx,1))=-F_{E_N}((1,-1-s_kx))$, the above equation 
determine the Eisenstein elements for the Eisenstein series $E_N$ completely. 
\end{proof}
From the above lemma, we conclude that $6F_N (g) = F_{E_N}(g)$. 
\begin{lemma}
\label{psi-values}
 For $E_N \in E_2(\Gamma_0(pq))$, let us consider the element $\sE_{E_N}^0$ of $ \HH_1(X_{\Gamma}-P_{-} , P_{+},\Z)$ defined by 
$\sE_{E_N}^0=\sum_{g \in \sP^1(\Z/pq\Z)} F_{E_N}(g) \xi^0(g)$. For all $c \in \HH_1(X_{\Gamma}-P_{+} , P_{-},\Z)$, we have
 $\sE_{E_N}^0 \circ c=\varphi_{E_N}(c)$
\end{lemma}
\begin{proof}
By Theorem~\ref{int-value}, we write the even Eisenstein element uniquely as
\[
\sum_{g \in \sP^1(\Z/pq\Z)} H_{E_N}(g) \xi^0(g). 
\]
By loc.\,cit.\,, $[g]_0 \circ [h]^0=1$ if and only if $\Gamma g=\Gamma h$. The function $H_{E_N}$ and $F_{E_N}$ coincide since  
$H_{E_N}(g)=\sum_{g \in \sP^1(\Z/pq\Z)} H_{E_N}(g) \xi^0(g) \circ \xi_0(g)= \sE^0_{E_N} \circ \xi_0(g)=F_{E_N}(g)$.
\end{proof}
For the modular curve $X_{\Gamma}$, we have a similar short exact sequence
\[
0  \rightarrow  \HH_1(X_{\Gamma}-P_{-} ,\Z) \rightarrow \HH_1(X_{\Gamma}-P_{-} , P_{+},\Z)  \xrightarrow{\delta^0} {\Z}^{P_{+}} \rightarrow  \Z \rightarrow 0.
\]
The boundary map $\delta^0$ takes a geodesic, joining the the point $r$ and $s$ of $P_{+}$ to the formal symbol $[r] - [s]$.

\section{Eisenstein elements and winding elements for \texorpdfstring{$\Gamma_0(pq)$}{X}}
\subsection{Eisenstein elements for \texorpdfstring{$\Gamma_0(pq)$}{X}}
We first prove an elementary number theoretic lemma. Recall, $l$ and $m$ are two unique integers such that
$lq + mp  \equiv 1 \pmod{pq}$ with $1 \leq l \leq (p - 1)$ and $1 \leq  m \leq (q - 1)$. 

\begin{lemma}
\label{twist}
For all $k$ with $1 \leq k \leq (q - 1)$, we can choose an integer $s(k) \in (\Z/q\Z)$ such that
\[
 (kp, -1) = (p, s(k)p -1)
\]
in $\sP^1(\Z/pq\Z)$. The map $k \rightarrow s(k)$ is a bijection $(\Z/q\Z)^* \rightarrow (\Z/q\Z) - \{\overline{m}\}$.
\end{lemma}
\begin{proof}
 For all $k$ with $1 \leq k \leq (q -1)$, let $k'$ be the inverse of $k$ in $(\Z/q\Z)^*$ . By Chinese remainder theorem,
we choose an unique $x$ with $1 \leq x \leq (pq - 1)$ such that $x \equiv -1 \pmod p$ and $x \equiv -k' \pmod q$. Observe
that $x$ is coprime to both $p$ and $q$. We write $x = s(k)p -1$ for an unique $s(k)$ with $0 \leq s(k) \leq (q -1)$. Since
$\Gamma_0(pq)\backslash \SL_2 (\Z) \cong \sP^1(\Z/pq\Z)$, we deduce that $(kp, -1) = (xkp, -x) = (-p, -x) = (p, x) = (p, s(k)p - 1)$
in $\sP^1(\Z/pq\Z)$.

 Consider the map $(\Z/q\Z)^* \rightarrow (\Z/q\Z)$ given by $k \rightarrow s(k)$. If $lq + mp  \equiv 1 \pmod{pq}$ then $m$ is not in the
image of this map. This map is one-one since if $s(k) = s(h)$ then $k \equiv h \pmod q$. Hence, the map
$(\Z/q\Z)^* \rightarrow (\Z/q\Z) - \{\overline{m}\}$ $k \rightarrow s(k)$ is a bijection.
\end{proof}
For all $t$ coprime to $pq$, consider the set $V$ of all matrices of the form $\alpha_t$. 
\begin{prop}
\label{boundary}
 The boundary of any element $$X =\sum_{g \in \sP^1(\Z/pq\Z)} F(g) [g]^*$$
 in $\HH_1(X_0(pq)-R \cup I, \partial(X_0(pq)),\Z)$ is of the form 
 \[
 \delta(X)=A(X)[\frac{1}{p}]+ B(X)[\frac{1}{q}]+C(X)[\infty]-(A(X)+B(X)+C(X)[0] 
\]
with
$$A(X)=\sum_{k=0}^{q-1} [F(\beta_k) - F(\beta_k S)], B(X) = \sum_{i=0}^{p-1} [F(\gamma_i) -F(\gamma_i S)]$$
and $C(X)=[F(0,1)-F(1,0)]$. 
\end{prop}
\begin{proof}
Choose an explicit coset representatives of $\Gamma_0(pq) \backslash \SL_2(\Z)$ (cf.  Lemma~\ref{explicit}) and write
$$X=C(X)[I]^*+\sum_{\alpha_t \in V} F(1,t) [\alpha_t]^*+\sum_{k=1}^{q-1}F(1, kp)[\alpha_{kp}]^*+\sum_{k=1}^{p-1} F(1,kq) [\alpha_{kq}]^*$$
$$+\sum_{i=0}^{q-1} F(p, ip-1) [\beta_i]^*+\sum_{j=0}^{p-1} F(q, jq-1) [\beta_j]^*.$$ 

According to Lemma~\ref{twist} for $1 \leq k \leq (q-1)$, we have $\alpha_{kp} S= Z \beta_{s(k)}$ for some $Z \in \Gamma_0(pq)$. We deduce that
\[
\sum_{k=1}^{q-1}F(1, kp)[\alpha_{kp}]^*+\sum_{i=0}^{q-1} F(p, ip-1) [\beta_i]^*=\sum_{k=1}^{q-1}(F(1,kp)[\alpha_{kp}]^*+
F(kp,-1)[\alpha_{kp} S]^*)+F(\beta_m)[\beta_m]^*
\]
and 
\[
 \sum_{k=1}^{p-1} F(1,kq) [\alpha_{kq}]^*+\sum_{j=0}^{p-1} F(q, jq-1) [\gamma_j]^*=\sum_{k=1}^{p-1} (F(1,kq) [\alpha_{kq}]^*+ F(kq,-1) [\alpha_{kq} S]^*)
 +F(\gamma_l)[\gamma_l]^*.
\]
A small check shows that $\delta([\alpha_{kp}]^*) = \delta([\alpha_p]^*)$ and $\delta([\alpha_{kp}]^*)= -\delta([\alpha_{kp} S]^*)$.

We now calculate $\delta([\beta_m]^*)$ and $\delta([\gamma_l]^*)$. 
Since $lq + mp \equiv 1 \pmod{pq}$ and $-I \in \Gamma_0(pq)$, so we get
\begin{equation}
 \left(\begin{smallmatrix}
1-q(l-1) & m(l-1)\\
(l-1)pq  & 1+lq(l-1)\\
\end{smallmatrix}\right)\left(\begin{smallmatrix}
m & -l\\
q  & p\\
\end{smallmatrix}\right)=\gamma \beta_m S
\end{equation}
 and 
\[
 \left(\begin{smallmatrix}
1-p(m+1) & -l(m+1)\\
(1+m)pq  & 1-mp(l+m)\\
\end{smallmatrix}\right)\left(\begin{smallmatrix}
m & -l\\
q  & p\\
\end{smallmatrix}\right)=\left(\begin{smallmatrix}
-1 & -l\\
q  & -mp\\
\end{smallmatrix}\right)=\gamma_l,
\]
for some $\gamma \in \Gamma_0(pq)$ and hence we have  $\Gamma_0(pq) \beta_m S=\gamma_l$. 
From $\delta([\beta_m]^*)= \delta([\alpha_q ]^*-[\alpha_p]^*)$
and $\delta([\gamma_l]^*) = \delta([\alpha_p]^* - [\alpha_q]^*)$, it is easy to see that
\[
\delta(\sum_{k=1}^{q-1} F(1,kp)[\alpha_{kp}]^*+ \sum_{i=0}^{q-1} F(p, jp-1)] [\beta_j]^*)
=\sum_{k=1}^{q-1} [F(1,kp)-F(kp,-1)] \delta([\alpha_p]^*)+ F(\beta_m)  \delta([\beta_m]^*).
\]
and 
\[
\delta(\sum_{k=1}^{p-1} F(1,kq)[\alpha_{kq}]^*+\sum_{j=0}^{p-1} F(q, jq-1)] [\gamma_j]^*)
=\sum_{k=1}^{p-1} [F(1,kq)-F(kq,-1)] \delta([\alpha_q]^*)+ F(q,lq-1)  \delta([\gamma_l]^*).
\]
\[
F(p,mp-1)  \delta([\beta_m]^*)+F(q,lq-1)  \delta([\gamma_l]^*)=[F(\beta_m)- F(\beta_mS)](\delta([\alpha_q]^*)-\delta([\alpha_p]^*)).
\]
Recall,  $\delta([\alpha_p]^*) =[0]-[\frac{1}{p}] $ and $\delta([\alpha_q]^*갱 ) =[0]-[\frac{1}{q}]$.
The above calculation shows that $\delta(X) = C(X)\delta([I]^* ) + A(X) \delta([\alpha_p]^* ) + B(X)\delta([\alpha_q]^*)$ with,
$A(X)=\sum_{k=0}^{q-1} [F(p,kp-1) - F(kp-1,-p)]$, $B(X) = \sum_{m=0}^{p-1} [F(\gamma'_l) -F(\gamma'_l S)]$
and $C(X)=F(I)-F(S)$.  We deduce the proposition.
\end{proof}
We also prove similar proposition for $\Gamma \subset \Gamma(2)$. 
\begin{prop}
\label{boundaryeven}
 The boundary of any element $$X =\sum_{g \in \sP^1(\Z/pq\Z)} F(g) \xi^0(g)$$
 in $\HH_1(X_\Gamma-P_{-}, P_{+},\Z)$ is of the form 
 \[
 \delta^0(X)=A'(X)[\frac{1}{p}]+ B'(X)[\frac{1}{q}]+C'(X)[\infty]-(A'(X)+B'(X)+C'(X)[0] 
\]
with
$$A'(X)=\sum_{k=0}^{q-1} [F(\beta'_k) -[\sum_{k=1}^{q-1} F(\alpha'_{kp})]-F(\gamma_l'),
B'(X) = \sum_{i=0}^{p-1} [F(\gamma'_i) -[\sum_{k=1}^{p-1} F(\alpha'_{kq})]-F(\beta_m')$$
and $C'(X)=[F(0,1)-F(\alpha'_{pq})]$. 
\end{prop}
\begin{proof}
 This is a straightforward calculation using the coset representatives of $\Gamma \backslash \Gamma(2)$ [cf.\, Lemma \ref{coset}].
\end{proof}
\begin{prop}
\label{boundact}
For $E \in \tE_{pq}$, the boundaries of almost Eisenstein elements $\sE'_E$ in $\HH_1(X_0(pq)-R \cup I, \partial(X_0(pq)),\Z)$ corresponding 
to the Eisenstein series $E$ are $-\delta(E)$ [\S  \ref{Eisensteinseries}].  
\end{prop}
\begin{proof}
For $E \in \tE_{pq}$, let $\sE'_E = \sum_{g \in \sP^1 (\Z/pq\Z)} G_E(g)[g]^*$ be the almost Eisenstein element. According
to Proposition~\ref{boundary}, we need to calculate $A(\sE'_E )$, $B(\sE'_E )$ and $C(\sE'_E )$.

For all $0 \leq k < (q - 1)$, $\beta_k T = \beta_{k+1}$ and $\beta_{q-1} T = \gamma \beta_0$ with 
$\gamma=\left(\begin{smallmatrix}
1+pq & q\\
-qp^2  & 1-qp\\
\end{smallmatrix}\right)$. We have an inclusion $\HH_1(Y_0(pq), \Z) \rightarrow \HH_1(Y_0(pq), R \cup I, \Z)$. Since 
$\{\rho^* , \gamma \rho^* \} =\{\beta_0\rho^* , \gamma \beta_0\rho^* \}=-\sum_{k=0}^{q-1} \{\beta_k \rho, \beta_k \rho^*\}$, 
we deduce that 
\[
\pi_E(\gamma)=\int_{z_0}^{\gamma z_0}E(z) dz=\sE'_E \circ \{z_0, \gamma z_0\}=-\sE'_E \circ (\sum_{k=0}^{q-1} \{\beta_k \rho, \beta_k \rho^*\})
=-\sum_{k=0}^{q-1} \sE'_E \circ \{\beta_k \rho, \beta_k \rho^*\}. 
\]
Applying Cor.~\ref{rho-int}, we have
$\sum_{k=0}^{q-1}\sE'_E \circ \{\beta_k \rho, \beta_k \rho^*갱\}=
\sum^{q-1}_{k=0} [G_E(\beta_k)-G_E(\beta_k S)] =-A(\sE'_E)$. 
Hence, we prove that $A(\sE'_E )=-\pi_E(\gamma)$. By interchanging $p$ and $q$, we have
$B(\sE'_E ) = -\pi_E(\gamma_0)$ for $\gamma_0 =
\left(\begin{smallmatrix}
1+pq & p\\
-p q^2  & 1-qp\\
\end{smallmatrix}\right).$

We now calculate  $\pi_E(\gamma)$ and $\pi_E(\gamma_0 )$ using \cite{MR800251}. Recall,
$\frac{1}{p}$ is a cusp with $e_{\Gamma_0(pq)}(\frac{1}{p})=q$. Consider 
the matrices $x=\left(\begin{smallmatrix}
1 & -q\\
-p  & 1+qp\\
\end{smallmatrix}\right)$ and $y=\left(\begin{smallmatrix}
1 & -p\\
-q  & 1+qp\\
\end{smallmatrix}\right)$ respectively. 
One can easily check that $x\left(\begin{smallmatrix}
1 & q\\
0  & 1\\
\end{smallmatrix}\right) x^{-1}= \gamma$ and $y\left(\begin{smallmatrix}
1 & p\\
0  & 1\\
\end{smallmatrix}\right) y^{-1}= \gamma_0$. Notice that $x(i \infty) = \Gamma_0(pq) \frac{1}{p}$ and $y(i \infty) = \Gamma_0(pq) \frac{1}{q}$.
By [\cite{MR800251}, p. 524], we deduce that $\pi_E(\gamma) = e_{\Gamma_0(pq)} ( \frac{1}{q}) a_0(E[\frac{1}{p}])$ and $\pi_{E_{pq}}(\gamma_0)= e_{\Gamma_0(pq)} (\frac{1}{p}) a_0(E[\frac{1}{p}])$.

According to Proposition~\ref{boundary}, the boundary of the almost Eisenstein element corresponding to an
Eisenstein series $E$ is
\[
\delta(\sE'_E)=A(\sE'_E)[\frac{1}{p}]+B(\sE'_E)[\frac{1}{q}]+C(\sE'_E)[\infty]-(A(\sE'_E)+B(\sE'_E)+C(\sE'_E))[0]
\]
with $A(\sE'_E)=qa_0(E[\frac{1}{p}])$, $B(\sE'_E)=pa_0(E[\frac{1}{q}])$ and $C(\sE'_E)=-[F(I)-F(S)]$. 
Applying Cor.~\ref{rho-int} again,  we deduce that
that $F(I)-F(S)=\int_{\rho}^{\rho^*}E(z)dz=-a_0(E)$. For $E \in E_2(\Gamma_0(pq))$, the boundary of 
$E$ is 
\[
\delta(E)=a_0(E)([\infty]-[0])+q a_0(E[\frac{1}{p}])([\frac{1}{p}]-[0])+pa_0(E[\frac{1}{q}])([\frac{1}{q}]-[0])
=\delta(\sE'_E).
\]
\end{proof}

Let $\beta$ and $h$ be the matrices $\left(\begin{smallmatrix}
1 & 2\\
0  & 1\\
\end{smallmatrix}\right)$ and $\left(\begin{smallmatrix}
1 & 1\\
0  & 2\\
\end{smallmatrix}\right)$
respectively. Let 
\[
 \pi_*:\HH_1(X_\Gamma-P_{-},P_{+},\Z) \rightarrow \HH_1(X_0(pq)-R \cup I, \partial(X_0(pq)),\Z)
\]
be the isomorphism defined by $\pi_*(\xi_0(g))=[g]^*$ [\cite{MR1363498}, Cor.\,1]. It is easy to see that $\delta(\pi_*(X))=\delta^0(X)$ for all 
$X \in \HH_1(X_\Gamma-P_{-},P_{+},\Z) $
\begin{prop}
\label{bound}
For all $E \in \tE_{pq}$, let $\sE_E^0 $ be the even Eisenstein element in $\HH_1(X_\Gamma-P_{-},P_{+},\Z)$ [\S ~\ref{Even}]. 
The boundary of the modular symbol $\pi_*(\sE^0_E)$ is $-6\delta(E)$.
\end{prop}
\begin{proof}
By Theorem 12, we explicitly write down the even Eisenstein element $\sE_E^0$ in the relative homology
group $\HH_1(X_{\Gamma} - P_{-}, P_{+}, \Z)$ as 
$$\sE_E^0 =\sum_{g \in \sP^1(\Z/pq\Z)} F_E(g)\xi_0 (g).$$
According to Proposition 19, we need to
calculate $A'( \sE_E^0), B'( \sE_E^0)$ and $C'( \sE_E^0)$. For $0 \leq k < (q -2)$, we have $\beta'_k \beta = \beta'_{k+2}$ . A small check shows
that $\beta'_{q-1} \beta = \beta'_1$ and $\beta'_{q-2} \beta = \gamma' \beta'_0$ with
\[
\gamma'=\left(\begin{smallmatrix}
1+2 pq(1+q) & 2q \\
-2q(p+pq)^2 & 1-2 pq(1+q)\\
\end{smallmatrix}\right) \in \Gamma.
\]
As a homology class in $\HH_1(X_{\Gamma} - P_{+}, P_{-} ,\Z)$, we have
\[
\{-1, \gamma'(-1)\} = \{\beta'_0 (-1), \gamma' \beta'_0(-1)\} = -\sum_{k=0}^{q-1}\{\beta'_k(1), \beta'_k (-1)\}=\sum_{k=0}^{q-1}\{\beta'_k(-1), \beta'_k (1)\}.
\]
By the definition of the even Eisenstein elements, we conclude that
\[
\int_{z_0}^{\gamma' z_0} k^*(\omega_E )= \sE_E^0 \circ \{z_0, \gamma' z_0 \} = - \sE_E^0\circ (\sum^{q-1}_{k=0}({\beta'_k(1), \beta'_k (-1)})
=-\sum^{q-1}_{k=0}\sE_E^0 \circ \{\beta'_k(1), \beta'_k (-1)\}.
\]
It is easy to see that $h AS B h^{-1} \in \SL_2(\Z)$ for all $A, B \in \Gamma(2)$.  Since $[\alpha'_{kq} S] = [\gamma'_{s(k)}]$ in $\sP^1(\Z/pq\Z)$, so
$\kappa'=\alpha'_{kq} S(\gamma'_{s(k)})^{-1} \in \Gamma_0(pq)$ and $h \kappa'h^{-1}\in \Gamma_0(pq)$. 
We deduce that the differential form
\[
k^*(\omega_E)=f(z)dz=[2 E(z)-\frac{1}{2} E(\frac{z+1}{2})] dz
\]
is invariant under $\kappa'$. According to the above argument,
\begin{equation}
\label{valuered}
F_E(\alpha'_{kq})=\int_{\alpha'_{kq}(1)}^{\alpha'_{kq}(-1)}f(z)dz=\int_{\alpha'_{kq}S(-1)}^{\alpha'_{kq}S(1)}f(z)dz=-\int_{\alpha'_{kq}S(1)}^{\alpha'_{kq}S(-1)}f(z)dz
\end{equation}
\[
=-\int_{\kappa'^{-1}\alpha'_{kq}S(1)}^{\kappa'^{-1}\alpha'_{kq}S(-1)}f(\kappa'z)d\kappa'z=-\int_{\gamma'_{s(k)}(1)}^{\gamma'_{s(k)}(-1)}f(z)dz=-F_E(\gamma'_{s(k)}).
\]
A similar calculation shows that $F_E (\gamma'_l ) = -F_E(\beta'_m)$ and $F_E(\alpha_{kp}) =-F_E( \beta_{s(k)})$ for some $s(k) \in (\Z/q\Z)^*$.
Applying Theorem~\ref{int-iso}, we have
\[
 \sum_{k=0}^{q-1} F_E(\beta'_k)= \sum_{k=0}^{q-1} \sE_E^0 \circ \{\beta_k'(1), \beta_k'(-1)\}=-\int_{z_0}^{\gamma' z_0} k^*(\omega_E).
\]
According to the definition of the period $\pi_E$ of the Eisenstein series $E(z)$ [cf. Section 3], we get
\[
 \int_{z_0}^{\gamma' z_0} k^*(\omega_E)=\int_{z_0}^{\gamma' z_0} [2 E(z)-\frac{1}{2} E(\frac{z+1}{2})] dz= 2 \pi_E(\gamma')-\pi_E(h\gamma'h^{-1}).
\]
We calculate $\pi_E(\gamma')$ and $\pi_E(h\gamma'h^{-1})$. From ~\ref{matrixmultiplication}, it is easy to see that $h\gamma'h^{-1}=\left(\begin{smallmatrix}
1+z & q v^2 \\
-4p^2q(1+q)^2 & 1-z\\
\end{smallmatrix}\right)$ with $v=(1- p(1 + q))$ and $z = 2pqv(1 + q)$. Furthermore, the matrix $h\gamma'h^{-1}$ decomposes as 
\[
h\gamma'h^{-1}=\left(\begin{smallmatrix}
1-p(1+q) & \frac{p(1+q)}{2} \\
-2p(1+q) & 1+p(1+q)\\
\end{smallmatrix}\right)
\left(\begin{smallmatrix}
1 & q  \\
0 & 1\\
\end{smallmatrix}\right)\left(\begin{smallmatrix}
1-p(1+q) & \frac{p(1+q)}{2} \\
-2p(1+q) & 1+p(1+q)\\
\end{smallmatrix}\right)^{-1}.
\] 
Since the matrix $\left(\begin{smallmatrix}
1-p(1+q) & \frac{p(1+q)}{2} \\
-2p(1+q) & 1+p(1+q)\\
\end{smallmatrix}\right)^{-1}$
takes the cusp $i \infty$ to $\frac{1}{p}$, we have $\pi_E(h \gamma'h^{-1})=q a_0(E[\frac{1}{p}])$. We further 
decompose $\gamma'$ as 
\[
\left(\begin{smallmatrix}
1 & -2q \\
-p(1+q) & 1+2pq(1+q)\\
\end{smallmatrix}\right)
\left(\begin{smallmatrix}
1 & 2q  \\
0 & 1\\
\end{smallmatrix}\right)\left(\begin{smallmatrix}
1 & -2q \\
-p(1+q) & 1+2pq(1+q)\\
\end{smallmatrix}\right)^{-1}.
\]
The matrix 
$\left(\begin{smallmatrix}
1 & -2q \\
p(1+q) & 1+2pq(1+q)\\
\end{smallmatrix}\right)$ takes the cusp $i \infty$ to $\frac{1}{p}$.  We deduce that $\pi_E(\gamma')=2 q a_0(E[\frac{1}{p}])$ and
$ \int_{z_0}^{\gamma' z_0} k^*(\omega_E)=3 a_0(E[\frac{1}{p}])$.  A simple  calculation shows that
\[
A'(\sE_E^0)=\sum_{k=0}^{q-1} F_E(\beta'_k)-\sum_{k=0}^{q-1} F_E(\alpha'_{kp})-F_E(\gamma'_m)=2\sum_{k=0}^{q-1} F_E(\beta'_k)=-6 a_0(E[\frac{1}{p}]).
\]
By interchanging $p$ and $q$, we obtain $B'(\sE_E^0)=-6 a_0(E[\frac{1}{q}])$. Since $\alpha'_{pq} S \in \Gamma_0(pq)$, a calculation similar to Equation~\ref{valuered} shows that
\[
F_E(I)=-F_E(\alpha_{pq})=\int_{1}^{-1} [2 E(z)-\frac{1}{2} E(\frac{z+1}{2})] dz=-\int_{-1}^{\beta(-1)} [2 E(z)-\frac{1}{2} E(\frac{z+1}{2})] dz=-3 a_0(E), 
\]
we conclude that $C'(\sE_E^0)=[F_E(I)-F_E(\alpha_{pq})]=-6 a_0(E)$ and hence $\delta^0(\sE_E^0)=\delta(\sE_E^0)=-6 \delta(E)$. 
\end{proof}
The inclusion map $i : (X_0(pq) - R \cup I, \partial(X_0(pq)) \rightarrow (X_0(pq), \partial(X_0 (pq))$ induces an onto map $i_* :
\HH_1(X_0(pq) - R \cup I, \partial(X_0(pq), \Z) \rightarrow \HH_1(X_0(pq), \partial(X_0(pq)), \Z)$ with $i_* ([g]^* ) = \xi(g)$. Note that $\delta([g]^*갱) =
[g.0] - [g.\infty] = \delta'(\xi(g)) =  \delta'(i_*([g]^*갱))$. From [\S~\ref{relative homology}], we have $\delta(c) = \delta'(i_*(c))$ for all homology class
$c \in  \HH_1 (X_0(pq)-R \cup I, \partial(X_0(pq), \Z)$.

\begin{lemma}
\label{integral}
 The integrals of every holomorphic differential on $X_0(pq)$ over $i_*(\sE_E')$  and $i_*\pi_*(\sE_E^0)$ are zero.
\end{lemma}
\begin{proof}
 A straightforward generalization of [\cite{MR1405312}, Lemma 5].
\end{proof}
We now prove the main Theorem ~\ref{Main-thm} of this article. 
\begin{proof}
By [\cite{MR1363498}, Cor.\,3], we obtain $i_*(\sE_E') \circ c = \sE_E' \circ i^* c = \int_c i_*(E(z)dz).$
Hence, $i_*(\sE'_E)$ is the Eisenstein element inside the space of modular symbols corresponding to $E$. 
By Proposition~\ref{boundact} and ~\ref{bound}, the boundary of $\pi_*(\sE^0_E)$ is same as the boundary of $6i_*(\sE'_E)$.

There is a non-degenerate bilinear pairing 
$ S_2(\Gamma_0(pq)) \times \HH_1(X_0(pq),\R) \rightarrow \C$
given by  $(f,c)=\int_c f(z) dz$. Hence, the integrals of the 
holomorphic differentials over $\HH_1(X_0(pq),\Z)$ are not always zero.
By Lemma~\ref{integral}, the integrals of every holomorphic differentials over $i_*(\sE'_E)$ and $i_*(\pi_*(\sE_E^0))$ are always zero. 
We deduce that 
\[
\sE_E=i_*(\sE'_E)=\frac{1}{6} i_*\pi_*(\sE_E^0)=\frac{1}{6} \sum_{g \in \sP^(\Z/pq\Z)} F_E(g) \xi(g).
\]
for  $E \in \tE_{pq}$.
Since $F_N(g) = \frac{1}{6} F_{E_N}(g)$, we obtain the theorem.
\end{proof}
\subsection{The winding elements of level $pq$}
\label{Winding}
Recall the concept of the {\it winding element}. 
\begin{definition}
\label{Winding}
[Winding elements]
Let $\{0,\infty\}$ denote the projection 
of the path from $0$ to $\infty$ in $\tH \cup {\sP}^1(\Q)$ to $X_0(pq)(\C)$. 
We have an isomorphism $\mathrm{H}_1(X_0(pq),\Z)\otimes \R = \mathrm{Hom}_{\C}(\mathrm{H}^0(X_0(pq),\Omega^1),\C)$.
Let $e_{pq} \in\mathrm{H}_1(X_0(pq),\R) $ corresponds to the homomorphism 
$\omega \rightarrow- \int_0^{\infty} \omega$. The modular symbol $e_{pq}$ is called the {\it winding element}.
\end{definition}
The winding elements are the elements of the space of modular symbols whose annihilators define ideals 
of the Hecke algebras with the $L$-functions of the corresponding quotients of the Jacobian non-zero.
In this paper, we found an explicit expression of the winding element. 
Let $e_{pq} \in  \HH_1 (X_0 (pq), \Z) \otimes \R$ be the winding element.
The following proposition help us to write down the winding element explicitly.
Since
$\displaystyle \sum_{x \in \partial(X_0(pq))}e_{\Gamma_0(pq))}(x)a_0(E[x])=0$, we write 
$$\delta(E)=  a_0(E) (\{ \infty \} - \{0 \}) + q a_0(E[\frac{1}{p}]) 
(\{ \frac{1}{p}\} - \{ 0 \}) + p a_0(E[\frac{1}{q}]) (\{ \frac{1}{q} \} - \{ 0 \}).$$
\begin{lemma}
\label{Epq}
The constant Fourier coefficients of $E_{pq}$ at cusps $0$, $\frac{1}{p}$, $\frac{1}{q}$ and $\infty$ are 
$\frac{1-pq}{24pq}$, $0$, $0$ and $\frac{pq-1}{24}$ respectively.
\end{lemma}
\begin{proof}
We first prove that the constant coefficient for the Fourier expansion of $E_{pq}$ at
the cusp $\frac{1}{p}$ is $0$. 
As usual, the constant term of the Fourier expansion of $E_{pq}$ at the cusp
$\frac{1}{p}$ is 
the constant term at $\infty$ of $E_{pq}[\beta_{0}]$. Similarly,  
the constant term of the Fourier expansion of $E_{pq}$ at the cusp
$\frac{1}{q}$ is 
the constant term at $\infty$ of $E_{pq}[\gamma_{0}]$.
Let $\Delta$ be the Ramanujan's cusp form of weight $12$. We write
$ \frac{d}{dz} \log \Delta(\beta(z))= 12 \frac{d}{dz} \log(pz+1)+ \frac{d}{dz}
\log\Delta(z)$ for $\beta$ = $\left(\begin{array}{cc}
1 & 0\\
p & 1\\
\end{array}\right)$.
A simple calculation shows that
$$\Delta(\frac{pq z}{pz +1})=\Delta(\left(\begin{array}{cc}
q & 0\\
1 & 1\\
\end{array}\right)pz)=\Delta(\left(\begin{array}{cc}
q & -1\\
1 & 0\\
\end{array}\right) \left(\begin{array}{cc}
1 & 1\\
0 & q\\
\end{array}\right)pz)=\\$$

$$\Delta(\left(\begin{array}{cc}
q & -1\\
1 & 0\\
\end{array}\right)(\frac{pz+1}{q}))= (\frac{pz+1}{q})^{12}\Delta(\frac{pz+1}{q}).$$
By taking logarithmic derivative, we deduce that 
$$\frac{d}{dz}log \Delta \left( \begin{array}{cc} q & -1\\ 1 & 0\\ \end{array}\right)
(\frac{pz+1}{q} )=12 \frac{d}{dz} \log(pz+1)+ \frac{d}{dz} \log \Delta(\frac{pz+1}{q}).$$
Since $E_{pq}(z)=\frac{1}{2\pi i}\frac{d}{dz} \log \frac{\Delta(pqz)} {\Delta(z)}$, the
above calculation 
shows that the constant term of $E_{pq}$ at the cusp $\frac{1}{p}$ is $0$.
Similarly, the constant term of $E_{pq}$ at the cusp $\frac{1}{q}$ is $0$.
The constant term of $E_{pq}$ is $\frac{pq-1}{24}$ at the cusp $\infty$ and
$\frac{1-pq}{24pq}$ at $0$.
\end{proof}
Using Lemma~\ref{integral} and Lemma~\ref{Epq},  we write
\begin{cor}
\label{cor}
 $$(1-pq) e_{pq}=\sum_{x \in(\Z/pq\Z)^*} F_{pq}((1,x)) \{0,\frac{1}{x}\}.$$
\end{cor}
\begin{remark}
For the Eisenstein series $E_p$ $\in$ $E_2(\Gamma_0(p))$, $\frac{1}{p}$ represents the cusp $\infty$ and
$\frac{1}{q}$ represents the cusp $0$. 
We deduce that $a_0(E_p[\beta_{0}])=\frac{p-1}{24}$ 
and $a_0(E_p[\gamma_{0}]) = \frac{1-p}{24p}$.
For the other Eisenstein series $E_q$ $\in$ $E_2(\Gamma_0(q))$, $\frac{1}{q}$ represents the cusp $\infty$ and
$\frac{1}{p}$ represents the cusp $0$. 
We deduce that $a_0(E_q[\gamma_{0}])=\frac{q-1}{24}$ and $a_0(E_q[\beta_{0}])=\frac{1-q}{24q}$.
\end{remark}
\bibliographystyle{crelle}
\bibliography{Eisensteinquestion.bib}
\end{document}